\documentclass[leqno,12pt]{article}
\usepackage{amsmath, amssymb}
\usepackage[dvips]{graphicx}
\usepackage{amsthm}
\title{ Special Liouville  metric with the Ricci condition  }
\author{ Katsuei  Kenmotsu }

\date{}

\newtheorem{thm}{Theorem}

\newtheorem{coro}{Corollary}
\newtheorem{remark}{Remark}
\begin{document}
\maketitle
\begin{abstract}
 Two necessary conditions for the induced metrics of parallel mean curvature surfaces in a complex space form of complex two-dimension are   observed.  One is similar to the Ricci condition of the classical surface theory in  Euclidean three-space and the other is related to the Liouville metric.   Conversely,  we prove that  a  special type of the Liouville  metric on a domain   in the Euclidean two-plane  whose   Gaussian curvature satisfies the  differential equation similar to the Ricci condition  is  explicitly determined by an elliptic function. We have  isometric  immersions from a simply connected two-dimensional Riemannian manifold  with the  special type of the  Liouville metric satisfying the Ricci condition to  the complex hyperbolic plane  with   parallel mean curvature vector.
\end{abstract}

\section{Introduction}
  The induced metric of a parallel mean curvature surface of a general type in the complex hyperbolic plane has the form of a special type of the Liouville metric (Kenmotsu \cite{ken2}).  In this paper,  we remark another necessary condition for  parallel mean curvature surfaces in complex space forms.  In  Section 2 of this paper, we prove  that  the Gaussian curvature of a parallel mean curvature surface in a two-dimensional complex space form  satisfies an elliptic differential equation that is intrinsic and  similar to the Ricci condition of  the classical surface theory (cf. Lawson \cite{laws}). 
  The main purpose of this study is to examine the converse of these two  facts.  Given  a domain in the  Euclidean two-plane $\mathbb{R}^{2}$,  we consider a special type  of the  Liouville   metric such that the Gaussian curvature satisfies the differential equation, which is the  same as that obtained in Section 2. Then, in Section 3,  we prove that  this  metric  is  explicitly determined by an elliptic function.  In the last section,   as an application of the result of Section 3,  we prove that a domain in $\mathbb{R}^{2}$ with the special type of the  Liouville metric such that the Gaussian curvature satisfies the  differential equation obtained in Section 2 can be isometrically immersed  in the complex hyperbolic plane  with parallel mean curvature vector. 
  
\section{Ricci condition on parallel mean curvature surfaces }
 In the 19th  century, Ricci observed  an intrinsic characterisation of constant mean curvature surfaces in $\mathbb{R}^{3}$, and in the late 20th century,  it was extended to   surfaces in  three-dimensional space forms by Lawson (Lawson \cite{laws}).  Later, this work was  generalised to parallel mean curvature surfaces in four-dimensional space forms by Eschenburg and Tribuzy \cite{esctri}, and Sakaki \cite{sakaki} and to minimal surfaces in   spheres with higher codimensions  by Vlachos \cite{vlacho}.

In this section, we demonstrate that the Gaussian curvature of  a parallel mean curvature surface in  a complex two-dimensional complex space form satisfies an elliptic   differential equation. This differential equation  is intrinsic and  similar to the Ricci condition for constant mean curvature surfaces in $\mathbb{R}^{3}$ (cf. Lawson \cite{laws}). For further studies on   the Ricci condition, we refer the reader to a recent paper by Tsouri and Vlachos \cite{tsou-vlach}.

Let $\overline{M}[4\rho]$ denote a complex  two-dimensional complex space form of  constant holomorphic sectional curvature $4\rho$,  and $M$ denote  an oriented and connected real two-dimensional Riemannian manifold.   Let $x\!:\!M\longrightarrow \overline{M}[4\rho]$ represent  an isometric immersion from $M$ to $\overline{M}[4\rho]$  with the Kaehler angle  $\alpha$ such that the  mean
curvature vector  $H$  is non-zero and parallel for the normal connection on the normal bundle
of the immersion.  Because the length of the mean curvature vector is constant,  
we may set  $|H| =2b >0$. This immersion is called   a  parallel mean curvature surface. The second fundamental form of $x$ is expressed by two complex-valued functions, $a$ and $c$ (Chern and Wolfson \cite{cherwolf}).  If  the Kaehler angle function is not constant and $a$ is not real-valued, then  the immersion $x$ is called of  a general type (Kenmotsu \cite{ken2}). Let $(u,v)$ represent a system of local coordinates in a neighbourhood  of  $M$ and  $f_{1}(u)\  (\mbox{resp}. \, f_{2}(v))$,  smooth real-valued functions of one variable $u\ (\mbox{resp}.\,  v)$,  which are considered as  functions on the neighbourhood, with  $f_{1}(u) +f_{2}(v) >0$ everywhere.  Then, the metric    $ds^{2}= (f_{1}(u) +f_{2}(v))(du^{2} +dv^{2})$ is called a  Liouville metric.  Refer studies by Kiyohara \cite{kiyo},  Knauf, Sinai, and Baladi \cite{ksb}  for the  results on Liouville metrics.   In particular, if $f_{2}(v) =0$ everywhere,  then   such a metric is called a special  Liouville metric.   Although the first fundamental form of any  surface of revolution in $\mathbb{R}^{3}$ is isometric to a special Liouville metric, we provide  new  insights into  special Liouville  metrics. Namely, we prove the following:
\begin{thm} Under the notation above, if $x\!: \!M\longrightarrow
\overline{M}[4\rho]$  is of a general type, then the Riemannian metric of $M$ is a  special Liouville metric,   the Gaussian curvature $K$ of  this metric satisfies  $ K  \leq 4b^{2}+2\rho $ and, at points where   $K <4b^{2}+2\rho $, 
\begin{equation}
\Delta\log\sqrt{4b^{2}+2\rho-K}=2K, 
\end{equation}
where $\Delta$ denotes the Laplace operator on $M$ with respect to the metric.
\end{thm}
Proof. We follow notation and formulas in \cite{ken2}.  By Lemma 2.4 of \cite{ken2}, there exists a complex-valued function
  $\mu(u)$ of a  real  variable $u$  such that 
$\phi = \mu(u)(du+ i dv) $,  where   $ds^{2} = \phi \bar{\phi} =  |\mu(u)|^{2}(du^{2} + dv^{2})$.  This proves the first part of Theorem 1.  
By (2.3) and (2.6) of \cite{ken2}, the Gauss equation of  $x$ is
\begin{equation}
|c|^{2} = \frac{1}{4}(4b^{2} + 2 \rho - K)  \ \geq 0,
\end{equation}
which proves the second part of Theorem 1.  We remark that if $x$ is of a general type, then $\alpha \neq 0$, and $ \alpha \neq \pi$, hence $\sin \alpha \neq 0$.
The Codazzi equation of $x$ is, 
by (2.5) of \cite{ken2}, 
$ dc\wedge\bar{\phi}=2c(a - b)\cot\alpha \, \phi\wedge\bar{\phi}.
$
\ At points where  $|c|>0$,  we have
$
\partial\log c =2(a-b) \cot\alpha \phi. 
$
Its  exterior differentiation  with   $(2.1),\ (2.2), \ (2.3), \ (2.4)$ of  \cite{ken2}  shows
\[ d(\partial\log c)=
\Big\{-\frac{K}{2}+2b(a-\bar{a})\frac{(1+\cos^{2}\alpha)}
{\sin^{2}\alpha} \Big\}\phi\wedge \bar{\phi}. 
\]
  Putting $c=|c|e^{\sqrt{-1}\eta}$, this  implies, at points where $|c| >0$,  
\begin{equation}
\Delta\log|c| = 2K,   
\end{equation}
where  the following formulas:  $4\partial \bar{\partial} \log |c| = \Delta \log |c| \phi \wedge \bar{\phi}$ and $\partial \bar{\partial}
 = -\bar{\partial} \partial$ are used.
The last part of Theorem 1  follows from $(2)$ and $(3)$, proving Theorem 1.
\vspace{0.5cm}

We provide three remarks about Theorem 1.

\begin{remark}

$(1)$ \  The first part of { \rm  Theorem 1}  holds under the assumption that $x$ is of a general type.  If  $x$ is not of a general type,  these parallel mean curvature surfaces are  classified as done by  {\rm Hirakawa \cite{hirakawa}}.

$(2)$\  It is proved by {\rm Kenmotsu \cite{ken3}}  that non-trivial parallel mean curvature surfaces of a general type  exist only when  $\rho = -3b^{2}$, so the Ricci condition $(1)$ is now,  $\Delta\log\sqrt{-2b^{2}-K}=2K,  \ (K < -2b^{2})$.

$(3)$\ The Ricci condition $(1)$  implies that the induced metric of a parallel mean curvature surface of a general type $x\!: \!M\longrightarrow
\overline{M}[4\rho]$ can be locally realized on a minimal surface in a $3$-dimensional real space form, $\mathbb{Q}^{3}_{\bar{c}}$, of curvature $\bar{c}=4b^{2} + 2 \rho$. This establishes a sort of Lawson correspondence between parallel mean curvature surfaces of a general type in $\overline{M}[4\rho]$ and minimal surfaces in $\mathbb{Q}^{3}_{\bar{c}}$.
\end{remark}

\section{Special Liouville metric with the Ricci condition  }
In this section, we study  the converse of Theorem 1.
Let $D$ denote an open, simply connected domain in $\mathbb{R}^{2}$.  In Theorem 2 and Corollary 1 of this section,  we prove that the Ricci condition  implies certain algebraic relation for a conformal factor and the Gaussian curvature  on $D$. Moreover, if the metric is a special Liouville metric, then Theorem 3 states that this is explicitly written by an elliptic function.

\begin{thm} Let  $ds^{2}$ denote a  Riemannian metric on $D$ and $b>0$ a positive number such that the Gaussian curvature of this metric  satisfies  $K < -2b^{2}$ and the Ricci condition  
\begin{equation}
\Delta\log\sqrt{-2b^{2}-K}=2K.
\end{equation}
  Then,  there exists a Riemannian metric $|\mu|^{2}|dw|^{2}$on $D$  which is isometric to $ds^{2}$, and satisfies 
\begin{equation}
 |\mu|^{2}\sqrt{-2b^{2}-K} =\mbox{constant} >0.
\end{equation}
\end{thm}
Proof.  The Euclidean  two-plane $\mathbb{R}^{2}$ is identified as the Gaussian plane $\mathbb{C}$,  so a point $(u,v)$ of $D$ is  written as  $z=u+iv \in \mathbb{C}, \ (i^{2}=-1)$. Let $z \in D$ denote an isothermal coordinate with $ds^{2}= |\nu|^{2}|dz|^{2}$, where $\nu$ is a non-zero  complex-valued function of $u$ and $v$.  By $2K= -\Delta \log |\nu|^{2}$ 
and (4), we have $\Delta \log( |\nu|^{2} \sqrt{-2b^{2}-K}) =0$. Since $ |\nu|^{2} \sqrt{-2b^{2}-K}$ is not zero on $D$, we can apply Lemma 3.12 of Eschenburg, Guadalupe and Tribuzy \cite{estri} to this function and therefore  there exists a holomorphic function $g(z)$ on $D$ such that $|g(z)| = |\nu|^{2}\sqrt{-2b^{2}-K} >0\ (z \in D)$. For any positive number  $c>0$,  let us define  a holomorphic transformation on $D$ by  $w=c^{-1/2} \int g^{1/2}dz$. Then we have $\nu dz =\mu dw$, where  $\mu = \nu c^{1/2} g^{-1/2}$ and  there exists a holomorphic function $k(w)$ such that
$ |k(w)| =|\mu| ^{2} \sqrt{-2b^{2}-K}$,  because $K$ is invariant by holomorphic transformations on $D$. By the definition of 
$\mu$, we have $|k(w)| =c  |g|^{-1} |g| = c $,  proving Theorem 2.

\begin{coro}
 In {\rm Theorem 2}, if $ds^{2}=|\nu|^{2}|dz|^{2}$ is a special Liouville metric, then so is $|\mu|^{2}|dw|^{2}$.
\end{coro}
Proof. We write $|\nu|^{2}= |\nu|^{2}(u)$. Then 
$ \Delta \log( |\nu|^{2}(u)\sqrt{-2b^{2}-K}) = 0 $ implies  $ d^{2}\log ( |\nu|^{2}(u)\sqrt{-2b^{2}-K})/du^{2}=0$,
because $K$ is a function of $u$ alone. So, we can write $ |\nu|^{2}(u) \sqrt{-2b^{2}-K} = c_{1}\exp(c_{2}u)$, for some $c_{1}\!>\!0$ and $c_{2} \in \mathbb{R}$. If $c_{2}=0$, then nothing to prove. When $c_{2}\neq 0$, we can take $g(z) = c_{1}\exp(c_{2}z/2)$ and $w= \int g^{1/2}dz$. Then, we have  $\mu =\nu c_{1}^{-1/2}\exp(-c_{2} z/4)$,
so  $|\mu|^{2} =| \nu|^{2} c_{1}^{-1}\exp(-c_{2}(z+\bar{z})/4) = |\nu|^{2}(u) c_{1}^{-1} \exp (-c_{2} u/2)$, which proves Corollary 1.

 Now we write the isothermal coordinate $w$ in  Corollary 1 as $w=u+iv, (i^{2}=-1)$, and put $ds^{2}=|\mu|^{2}|dw|^{2} =\lambda^{2}(u)(du^{2}+dv^{2}), \ (\lambda(u) >0, (u,v) \in D)$.    By (5),  we have 
 $\lambda^{4}(u)(-2b^{2}-K) = c_{1} , \ (c_{1} >0)$.  It follows from $K=(-\lambda \lambda'' + \lambda'^{2})/\lambda^{4}$ that
$$ 
\lambda''(u) \lambda(u) - \lambda'(u)^{2}  = c_{1} + 2b^{2} \lambda(u)^{4}, \ (c_{1} >0).
$$ 
It is integrated as
\begin{equation}
\lambda'(u)^{2} = -c_{1} + c_{2}\lambda(u)^2 + 2b^{2}\lambda(u)^{4}, \ (c_{2} \in R),
\end{equation}
where $c_{2}$ is an integral constant.
The right hand side of (6)   is decomposed as 
 $\lambda'(u)^{2} = 2b^{2}(\lambda^{2} -\lambda_{-})(\lambda^{2} - \lambda_{+})$, where
 $\lambda_{+}$ and $\lambda_{-}$ are real numbers such that  $\lambda_{+} >0> \lambda_{-}$ because of  $c_{1} >0$. Hence it holds that
   $\lambda^{2} \geq  \lambda_{+}>0$.  Let us define $\theta =\theta(u)$ by
\begin{equation}
\cos \theta =\frac{ \sqrt{\lambda_{+}}}{\lambda(u)},
\end{equation}
which  yields, by (6), 
$$ 
\theta'(u)^{2} = 2b^{2}(\lambda_{+} - \lambda_{-})(1+ \frac{\lambda_{-}}{\lambda_{+}-\lambda_{-}} \sin^{2}\theta).
$$ 
In other words,
\begin{equation}
\theta'(u)^{2} = \sqrt{c_{2}^{2}+ 8b^{2}c_{1}}- \left(\frac{c_{2}+  \sqrt{c_{2}^{2}+ 8b^{2}c_{1}}}{2} \right)\sin^{2}\theta .
\end{equation}
Now we introduce the constant $k^{2}$ and two functions $\tilde{u},\   \tilde{\theta}(\tilde{u})$.
$$ 
 k^{2}= \frac{c_{2}+\sqrt{c_{2}^{2}+ 8b^{2}c_{1}}}{2\sqrt{c_{2}^{2}+ 8b^{2}c_{1}}}, \ \tilde{u}= (c_{2}^{2}+  8b^{2}c_{1})^{\frac{1}{4}}u,     
    \  \tilde{\theta}(\tilde{u})= \theta( (c_{2}^{2}+  8b^{2}c_{1})^{-\frac{1}{4}}\tilde{u}).
$$ 
Then, we have  $(d\tilde{\theta}/d\tilde{u})^{2} = 1 - k^{2} \sin^{2}\tilde{\theta}$, that is,
\begin{equation}
\theta(u) = am ((c_{2}^{2} +  8b^{2}c_{1})^{\frac{1}{4}}u,k),
\end{equation}
where $am(\cdot,k)$ is  Jacobi's amplitude with the modulus $k$.  Hence, we have, by (7) and (9),
\begin{equation}
\lambda(u)= \frac{ \sqrt{\lambda_{+}}}{\cos\theta(u)}=\frac{ \sqrt{\lambda_{+}}}{cn((c_{2}^{2}+8b^{2}c_{1})^{1/4} u,k)} , 
\end{equation}
where $4b^2\lambda_{+}= -c_{2}+\sqrt{c_{2}^{2}+8b^{2}c_{1}}$ and  $cn(\cdot,k)$ is the Jacobi elliptic function with the modulus $k$.
Thus, we  have the two-parameter family of  Jacobi elliptic functions. The main result of this study is the following:

\begin{thm} Let $D$ denote an open, simply connected domain in $\mathbb{R}^{2}$ and $ds^{2}$ a special Liouville metric on $D$.  Suppose that the Gaussian curvature of this  metric   satisfies  $K < -2b^{2}$ everywhere  for some $b>0$   and the  Ricci  condition $(4)$. Then,  there exists a special Liouville metric
 on $D$  which is isometric to $ds^{2}$,  and  whose conformal factor  is  explicitly determined by a Jacobi elliptic function and two  real constants. 
\end{thm}

\section{Application to parallel mean curvature surfaces}
In this section, we apply  the results of the previous section to obtain immersions from a domain $D$ in
 $\mathbb{R}^{2}$ 
to a complex two-dimensional complex space form. For simplicity, we assume that $6b^{2}=1$.  Let $\mathbb{CH}^{2}[-2]$ be the complex hyperbolic plane of   constant holomorphic sectional curvature $-2$. 
\begin{thm}
Let $D$ denote an open, simply  connected domain in $\mathbb{R}^{2}$,  $ds^{2}$ a special Liouville metric on $D$ so that the Gaussian curvature of this metric satisfies $K < -1/3$  and $\Delta\log\sqrt{-1/3-K}=2K$. Then, there exists a family of  isometric immersions from $(D,  ds^{2})$ into $\mathbb{CH}^{2}[-2]$ such that the mean curvature vector  $H$  of the immersion in the family  is parallel for the normal connection of the normal bundle of the immersion and $|H| = 2/\sqrt{6}$.
\end{thm}
Proof.   We  find certain one-parameter sub-family of   Jacobi elliptic functions obtained in Theorem 3. For a given $c_{1}>0$, put  $c_{2}=c_{1}/6-2$. Then, the formula  (9) becomes  $\theta'(u)^{2} = 2+c_{1}/6 - c_{1}/6  \sin^{2}\theta$. Therefore when $0<c_{1}<3/2$, setting $p=c_{1}/6$ and $\theta = \gamma$, we have the equation (4.3) of \cite{ken3}.  Take $\alpha$ such that $3\cos \alpha= -\sin \gamma $.     Using this $\alpha$ and putting $b=1/\sqrt{6}$ let us define $\mu, a,$ and $c$ by (4.8) of \cite{ken3}. Then the metric defined by $\mu$ is isometric to $ds^{2}$ and by Theorem 4.1 of \cite{ken3},  Theorem 4 is proved when  $0<c_{1}<3/2$.  When  $3/2 <c_{1} $, setting $c_{2}= 2-c_{1}/6, \ p=c_{1}/6$ and $\theta =\tilde{\gamma}$, we have the equation (4.5) of \cite{ken3}.  In  the same way as above, we proved Theorem 4.

\section{Added in proof}
\leftline{  {\bf Surface of revolution in $\mathbb{R}^{3}$}} \quad  The first fundamental form of a surface of revolution in $\mathbb{R}^{3}$ is isometric to a special Liouville metric. In fact, let $I = ds^{2} + y(s)^{2}dv^{2}, \ (y(s)>0)$ denote the first fundamental form of a surface of revolution
 $(x(s),y(s)\cos v,y(s)\sin v) \subset \mathbb{R}^{3}$, where $s$ is the arc length parameter of the profile curve. Then we have $I=\lambda(u)^{2}(du^{2}+dv^{2})$, where $u=u(s) := \int ds/y(s)$, $s=s(u)$ is its inverse function, and $\lambda(u) := y(s(u))$. We remark that the conformal factor $\lambda(u)$ must satisfy $\lambda^{2}(u) \geq \lambda'(u)^{2}$  everywhere because of $y'(s)^{2}\leq 1$. Conversely the special Liouville metric $I=\lambda(u)^{2}(du^{2}+dv^{2})$ with  $\lambda^{2}(u) \geq \lambda'(u)^{2}$ induces naturally a surface of revolution in $R^{3}$ where the profile curve  is defined by $y(u):= \lambda(u),\ x(u): = \int \sqrt{\lambda^{2}(u) - \lambda'^{2}(u)}du$.  In particular, if the special Liouville metric satisfies the Ricci condition, then $\lambda(u)$ is determined by  (10).  The shape of the surface resulting from the profile curve resembles the front part  of a trumpet,  because the both coordinate functions of the profile curve  are monotone increasing  on an interval of $u>0$.
\vspace{0.3cm}

\leftline{{\bf Open problem}}  The results by Eschenburg and Tribuzy \cite{esctri} indicate a problem of determining  whether the hypothesis in Theorem 4, which is the special Liouville  metric, is necessary.

\medskip

\begin{flushleft}

 Katsuei  Kenmotsu \\
Mathematical Institute,  Tohoku University  \\
980-8578 \quad  Sendai, Japan \\
email:  kenmotsu-math@tohoku.ac.jp
\end{flushleft}
\end{document}